\documentclass[10pt]{gsm-l}

\usepackage{amssymb}

\usepackage{graphicx}

\usepackage{amsmath}
\usepackage{amssymb}

\usepackage[alphabetic,initials,lite,nobysame]{amsrefs}
\usepackage{url}

\makeatletter
\newdimen\@defaultposttitleskip  \@defaultposttitleskip=7.5pc
\newdimen\@posttitleskip
\newif\if@epigraph
\newtoks\@epigraph
\long\def\epigraph#1{\@epigraphtrue
	\@epigraph=\@xp{\the\@epigraph\endgraf\vspace{.5\baselineskip}\itshape#1}}
\long\def\@makechapterhead#1{\global\topskip\normaltopskip
	\begingroup
	\parskip\z@skip
	\vbox to\topskip{%
		\chapter@number
		\vss
	}\penalty\@M
	\Huge\bfseries \hsize24pc \raggedright
	\top@space{11pc}%
	\noindent\ignorespaces #1\par \endgroup
	\@posttitleskip=\@defaultposttitleskip
	\if@epigraph
	\@posttitleskip=3.5pc
	\vspace{1.5pc}
	\begingroup \small
	\rightskip1pc\relax \leftskip.20\hsize\relax
	\parfillskip-1pc\relax \parindent\z@\relax
	\def\\{\hfill\vadjust{\vskip1pt}\break
		\null\hfill\rm---\thinspace\ignorespaces}%
	\itshape\the\@epigraph\endgraf
	\endgroup
	\global\@epigraph{}\global\@epigraphfalse
	\fi
	\if@index
	\ifx\@empty\indexintro
	\bb@space{\@posttitleskip}%
	\else
	\vspace{3pc}%
	\begingroup \small
	\parbox[t]{27pc}{\leftskip3pc\normalfont\indexintro\par}%
	\endgroup
	\bb@space{3pc}%
	\fi
	\else
	\@dropfolio
	\bb@space{\@posttitleskip}%
	\fi
	\@afterheading
} \makeatother

\newtheorem{theorem}{Theorem}[chapter]

\theoremstyle{definition}

\theoremstyle{remark}

\numberwithin{section}{chapter}
\numberwithin{equation}{chapter}

\makeindex

\begin{document}

\frontmatter

\title[Fractional Sobolev Spaces]{A First Course in Fractional Sobolev Spaces}


\author{Giovanni Leoni}
\address{Department of Mathematical Sciences, Carnegie Mellon University, Pittsburgh,
	PA, 15213}
\email{giovanni@andrew.cmu.edu}
\urladdr{http://www.math.cmu.edu/people/fac/leoni.html}
\thanks{NSF}



\maketitle

\cleardoublepage
\thispagestyle{empty}
\vspace*{13.5pc}
\begin{center}
 Dedicated to Irene, with friendship, \\[2pt] and to the memory of all the victims of
 the pandemic
\end{center}
\cleardoublepage

\setcounter{page}{7}

\chapter*{Contents}

\noindent Preface \hfill xi

\noindent \textbf{Part 1. Fractional Sobolev Spaces in One Dimension}

\noindent Chapter 1. Fractional Sobolev Spaces in One Dimension \hfill 3

\S 1.1. Some Useful Inequalities and Identities \hfill 4

\S 1.2. Fractional Sobolev Spaces \hfill 12

\S 1.3. Extensions \hfill 22

\S 1.4. Some Equivalent Seminorms \hfill 27

\S 1.5. Density \hfill 33

\S 1.6. The Space $W_{0}^{s,p}(I)$ \hfill 38

\S 1.7. Hardy's Inequality \hfill 41

\S 1.8. The Space $W_{00}^{s,p}(I)$ \hfill 46

\S 1.9. Notes \hfill 48

\noindent Chapter 2. Embeddings and Interpolation \hfill 49

\S 2.1. Embeddings: The Endpoints \hfill 50

\S 2.2. Embeddings: The General Case \hfill 64

\S 2.3. Interpolation Inequalities \hfill 74

\S 2.4. Notes \hfill 86

\noindent Chapter 3. A Bit of Wavelets \hfill 87

\S 3.1. Weighted $\ell^{p}$ Spaces \hfill 88

\S 3.2. Wavelets in $L^{p}$ \hfill 89

\S 3.3. Wavelets in $L^{1}$ \hfill 92

\pagebreak

\S 3.4. Wavelets in $W^{1,1}$ \hfill 93

\S 3.5. Wavelets in $W^{s,p}$ \hfill 96

\S 3.6. Back to Interpolation Inequalities \hfill 101

\S 3.7. Notes \hfill 107

\noindent Chapter 4. Rearrangements \hfill 109

\S 4.1. Polarization \hfill 109

\S 4.2. Polarization in $W^{s,p}(\mathbb{R})$ \hfill 124

\S 4.3. Symmetric Decreasing Rearrangement \hfill 127

\S 4.4. Symmetric Decreasing Rearrangement in $W^{s,p}(\mathbb{R})$ \hfill 142

\S 4.5. Notes \hfill 144

\noindent Chapter 5. Higher Order Fractional Sobolev Spaces in One Dimension \hfill 145

\S 5.1. Higher Order Fractional Sobolev Spaces \hfill 145

\S 5.2. Extensions \hfill 151

\S 5.3. Density \hfill 158

\S 5.4. Some Equivalent Seminorms \hfill 159

\S 5.5. Hardy's Inequality \hfill 165

\S 5.6. Truncation \hfill 168

\S 5.7. Embeddings \hfill 170

\S 5.8. Interpolation Inequalities \hfill 174

\S 5.9. Notes \hfill 181

\noindent \textbf{Part 2. Fractional Sobolev Spaces}

\noindent Chapter 6. Fractional Sobolev Spaces \hfill 185

\S 6.1. Definition and Main Properties \hfill 185

\S 6.2. Slicing \hfill 201

\S 6.3. Some Equivalent Seminorms \hfill 214

\S 6.4. Density of Smooth Functions \hfill 222

\S 6.5. The Space $W_{0}^{s,p}$ and Its Dual \hfill 228

\S 6.6. Hardy's Inequality \hfill 234

\S 6.7. The Space $W_{00}^{s,p}$ \hfill 244

\S 6.8. Density in $\dot{W}^{s,p}$ \hfill 251

\S 6.9. Notes \hfill 255

\noindent Chapter 7. Embeddings and Interpolation \hfill 257

\S 7.1. Campanato Spaces \hfill 258

\pagebreak

\S 7.2. Embeddings: The Subcritical Case \hfill 261

\S 7.3. Embeddings: The Critical Case \hfill 264

\S 7.4. Embeddings: The Supercritical Case \hfill 275

\S 7.5. Embeddings: The General Case \hfill 277

\S 7.6. The Limit of $W^{s,p}$ as $s\rightarrow1^-$ \hfill 296

\S 7.7. Interpolation Inequalities \hfill 300

\S 7.8. Notes \hfill 308

\noindent Chapter 8. Further Properties \hfill 309

\S 8.1. Extension: Lipschitz Domains \hfill 309

\S 8.2. Extension: The General Case \hfill 313

\S 8.3. Derivatives \hfill 329

\S 8.4. Embeddings and Interpolation Inequalities \hfill 331

\S 8.5. Notes \hfill 335

\noindent Chapter 9. Trace Theory \hfill 337

\S 9.1. Traces of Weighted Sobolev Spaces \hfill 337

\S 9.2. The Trace Operator in $W^{s,p}$ \hfill 349

\S 9.3. Half-Spaces \hfill 355

\S 9.4. Special Lipschitz Domains \hfill 361

\S 9.5. Bounded Lipschitz Domains \hfill 366

\S 9.6. Unbounded Lipschitz Domains \hfill 372

\S 9.7. Notes \hfill 384

\noindent Chapter 10. Symmetrization \hfill 387

\S 10.1. Polarization \hfill 387

\S 10.2. Polarization in $W^{s,p}(\mathbb{R}^{N})$ \hfill 391

\S 10.3. Spherically Symmetric Rearrangement \hfill 394

\S 10.4. Spherical Symmetric Rearrangement in $W^{s,p}(\mathbb{R}^{N})$ \hfill 402

\S 10.5. Notes \hfill 405

\noindent Chapter 11. Higher Order Fractional Sobolev Spaces \hfill 407

\S 11.1. Definition \hfill 407

\S 11.2. Some Equivalent Seminorms \hfill 415

\S 11.3. Slicing \hfill 422

\S 11.4. Embeddings \hfill 426

\S 11.5. Interpolation Inequalities \hfill 436

\S 11.6. Superposition \hfill 448

\pagebreak

\S 11.7. Extension \hfill 451

\S 11.8. Trace Theory \hfill 454

\S 11.9. Notes \hfill 479

\noindent Chapter 12. Some Equivalent Seminorms \hfill 481

\S 12.1. The Interpolation Seminorm \hfill 481

\S 12.2. The Littlewood--Paley Seminorm \hfill 490

\S 12.3. Notes \hfill 518

\noindent \textbf{Part 3. Applications}

\noindent Chapter 13. Interior Regularity for the Poisson Problem \hfill 521

\S 13.1. Cacciopoli's Inequality \hfill 521

\S 13.2. $H^{2}$ Regularity \hfill 525

\S 13.3. $W^{2,p}$ Regularity \hfill 532

\S 13.4. $W^{s,p}$ Regularity \hfill 541

\S 13.5. Notes \hfill 545

\noindent Chapter 14. The Fractional Laplacian \hfill 547

\S 14.1. Definition and Main Properties \hfill 547

\S 14.2. The Riesz Potential \hfill 556

\S 14.3. Existence of Weak Solutions \hfill 561

\S 14.4. An Extension Result \hfill 567

\S 14.5. Notes \hfill 575

\noindent Bibliography \hfill 577

\noindent Index \hfill 585


\chapter*{Preface}

What are fractional Sobolev spaces and why are they important? In his blog,
Terence Tao gives a beautiful answer to these questions: If you consider a
simple function $u$ of the form $u=h\chi_{E}$, where $E\subseteq\mathbb{R}%
^{N}$ is a Lebesgue measurable set and $h$ is a positive constant, then the
$L^{p}$ norm of $u$ is given by the number $h(\mathcal{L}^{N}(E))^{1/p}$,
where $\mathcal{L}^{N}(E)$ is the volume of the set $E$. Hence, you measure
the height $h$ of the function $u$ and its width. Sobolev spaces also measure
the regularity of the function and how quickly it oscillates (the frequency
scale of $u$).

A function $u$ belongs to the Sobolev space $W^{1,p}(\mathbb{R}^{N})$ if $u$
belongs to $L^{p}(\mathbb{R}^{N})$, and its (weak) gradient $\nabla u$ belongs
to $L^{p}(\mathbb{R}^{N};\mathbb{R}^{N})$. When $p=2$, it is possible to give
an equivalent definition using Fourier transforms. Indeed, you can check that
$u$ belongs to $W^{1,2}(\mathbb{R}^{N})$ if and only if its Fourier transform
$\hat{u}$ and $\Vert\xi\Vert\hat{u}$ belong to $L^{2}(\mathbb{R}^{N})$, with%
\begin{align*}
\int_{\mathbb{R}^{N}}|u(x)|^{2}dx  &  =\int_{\mathbb{R}^{N}}|\hat{u}(\xi
)|^{2}d\xi,\\
\int_{\mathbb{R}^{N}}\Vert\nabla u(x)\Vert^{2}dx  &  =\int_{\mathbb{R}^{N}%
}(2\pi\Vert\xi\Vert)^{2}|\hat{u}(\xi)|^{2}d\xi.
\end{align*}
To define the fractional Sobolev space $W^{s,2}(\mathbb{R}^{N})$, replace
$\Vert\xi\Vert$ with $\Vert\xi\Vert^{s}$, where $0<s<1$. Intuitively, this is
saying that instead of having a full derivative, you only have a fraction of a
derivative. After some computations, you find the following relations%
\begin{align*}
\int_{\mathbb{R}^{N}}|u(x)|^{2}dx  &  =\int_{\mathbb{R}^{N}}|\hat{u}(\xi
)|^{2}d\xi,\\
\int_{\mathbb{R}^{N}}\int_{\mathbb{R}^{N}}\frac{|u(x)-u(y)|^{2}}{\Vert
x-y\Vert^{N+s2}}dx  &  =C\int_{\mathbb{R}^{N}}\Vert\xi\Vert^{2s}|\hat{u}%
(\xi)|^{2}d\xi,
\end{align*}
where $C=C(N,s)>0$. You then extend this definition to the case $p\neq2$. You
say that a function $u\in L^{p}(\mathbb{R}^{N})$ belongs to the fractional
Sobolev space $W^{s,p}(\mathbb{R}^{N})$ if
\begin{equation}
\int_{\mathbb{R}^{N}}\int_{\mathbb{R}^{N}}\frac{|u(x)-u(y)|^{p}}{\Vert
x-y\Vert^{N+sp}}dx<\infty. \label{seminorm}%
\end{equation}
All this is quite neat, but it still does not answer why you should care. It
turns out that fractional Sobolev spaces play a central role in the calculus
of variations, partial differential equations, and harmonic analysis. Let me
try to explain.

In 1957 Gagliardo \cite{gagliardo1957} proved that, when $1<p<\infty$, the
fractional Sobolev space $W^{1-1/p,p}(\mathbb{R}^{N-1})$ is the trace space of
the Sobolev space $W^{1,p}(\mathbb{R}_{+}^{N})$, with $\mathbb{R}_{+}^{N}$ the
open half-space $\{x=(x_{1},\ldots,x_{N}):\,x_{N}>0\}$. The case $p=2$ had
been studied by Aronszajn \cite{aronszajn1955}, Prodi \cite{prodi1956}, and
Slobodecki\u{\i} \cite{slobodeckii1958}. For a regular function, the trace is
the restriction of the function to the boundary. This result is significant
because if you want to establish the existence of a weak solution to partial
differential equations with non homogeneous boundary conditions, you need to
consider fractional Sobolev spaces. Important examples are the non homogeneous
Dirichlet or Neumann problem for the Poisson equation or the heat equation.

The same happens if you want to find a minimizer or a saddle point of a
functional defined on a class of functions having the same inhomogeneous trace
on the boundary. This problem is central in the calculus of variations.

Besides existence, fractional Sobolev spaces are essential when studying the
regularity of solutions of partial differential equations. Imagine that you
proved a differential function satisfies the differential equation $u^{\prime
}=f(u)$, where $f$ is a smooth function. But then the right-hand side $f(u)$
is a differentiable function (by the chain rule), and so you learn that
$u^{\prime}$ is differentiable itself. You have gained one extra derivative.
In turn, $f(u)$ is twice differentiable and, using the equation, $u^{\prime}$
becomes twice differentiable. You can continue in this way, depending on how
smooth the function $f$ is. This procedure is a \emph{bootstrap method}. Now,
if you replace the right-hand side $f(u)$ with something more complicated, or
if instead of an ordinary differential equation, you have a partial
differential equation, then you might not be able to gain a full derivative.
But in some cases, you might be able to get a fraction of a derivative. So
again, this brings you to fractional Sobolev spaces.

To give an example, the existence and smoothness of solutions to
Navier--Stokes equations in $\mathbb{R}^{3}$ are among the Millenium Problems
of the Clay Institute. While we are still far from finding a solution,
fractional Sobolev spaces have played a significant role in the theory of the
Navier--Stokes equations \cite[Chapter 5]{bahouri-chemin-danchin-book2011},
\cite[Theorem 3.5]{chemin-desjardins-gallagher-grenier-book2006}, \cite[Part
5]{lemarie-book2002}, \cite{sohr-book2001}.

Finally, fractional Sobolev spaces and Besov spaces are ubiquitous in harmonic
analysis and singular integrals. Indeed, it is possible to characterize these
spaces using the Littlewood--Paley decomposition. In a way, Besov cases
$B_{q}^{s,p}$ are the grown-up version of fractional Sobolev spaces since they
have an extra parameter $q$. When $0<s<1$ and $p=q$, $B_{p}^{s,p}%
(\mathbb{R}^{N})=W^{s,p}(\mathbb{R}^{N})$. Hence, fractional Sobolev spaces
are a subclass of Besov spaces.

And so now we come to the sticky point. Why do we need another book on
fractional Sobolev spaces? There are several ways to characterize Besov spaces
and fractional Sobolev spaces. Besides the intrinsic formulation given in
(\ref{seminorm}), you can use abstract interpolation theory to represent
$W^{s,p}(\mathbb{R}^{N})$ as an intermediate space between $L^{p}%
(\mathbb{R}^{N})$ and $W^{1,p}(\mathbb{R}^{N})$. This method has the advantage
that abstract interpolation theory lets you deduce several theorems (for
example, embedding theorems) from the corresponding ones for $W^{1,p}%
(\mathbb{R}^{N})$. And then there is also the Littlewood--Paley decomposition,
which carries with it all the power of Fourier analysis (for example, Fourier
multipliers). Leading experts navigate seamlessly from one formulation to the
next to get the most elegant proof. In a way, this is great. I don't know of
any other area where such different disciplines all converge.

But at the same time, if you are an advanced undergraduate student or a
beginning graduate student, reading Russian literature
\cite{besov-ilin-nikolski-book1979}\ or Triebel's monographs
\cite{triebel-book1995}, \cite{triebel-book2010}, \cite{triebel-book1992},
\cite{triebel-book2006}\ can be challenging. I wrote this book as a gentle
introduction to fractional Sobolev spaces for advanced undergraduate students
(Part I) and graduate students (Part II). In a way, my book is like
\textit{The Hobbit} and \cite{besov-ilin-nikolski-book1979},
\cite{triebel-book1995}, \cite{triebel-book2010}, \cite{triebel-book1992},
\cite{triebel-book2006} are \textit{The Lord of the Rings}.

With the exception of Chapter 12, all the proofs in this
book rely on the intrinsic seminorm (\ref{seminorm}). Since the
\textquotedblleft intrinsic" approach bypasses the need for harmonic analysis,
it's more readily amenable to the situations we encounter in studying boundary
value problems. To my knowledge, no one has assembled the intrinsic techniques
until now.

As for why this book now? In 2007 Caffarelli and Silvestre published their
paper \cite{caffarelli-silvestre2007}, in which they studied properties of
solutions to the fractional Laplacian. Their article, which has more than two
thousand citations in Google Scholar, started a ripple in mathematics.
Mathematicians extended several results from the theory of elliptic partial
differential equations (the theory of viscosity solutions, critical point
theory, etc.) and the calculus of variations to nonlocal operators.

I should make clear that this book is not on the fractional Laplacian. I wrote
one chapter about it because if you consider the Euler--Lagrange equations of
the fractional seminorm (\ref{seminorm}), you do get the fractional Laplacian.
But this is as far as I went.

As a byproduct of this new interest on the fractional Laplacian, new proofs of
classical results on fractional Sobolev spaces appeared in the recent literature.

When the pandemic hit (in March 2020), I suddenly found myself stuck at home
with two cats (Gauss and Hilbert). While they kept me sane during the months
of lockdown, I was unable to concentrate on my research. Try focusing when you
have two cats squabbling (for hours!) on who gets to sleep on your lap. So
instead, I devoted myself to online shopping (alas) and navigating the
internet. Early enough, I stumbled on the paper of Nguyen, Diaz, and Nguyen
\cite{nguyen-diaz-nguyen2020}, in which they gave a new proof of a classical
embedding for fractional Sobolev spaces using only maximal functions. When I
was writing the chapter on Besov spaces for the first edition of my book
\cite{leoni-book2009}, I had spent months trying to find in the literature
precisely these types of proofs. For the most part, I was not satisfied with
the final result. In the second edition \cite{leoni-book2017}, I changed that
chapter completely. I caved in and used abstract interpolation and the
Littlewood--Paley theory. But I was still not happy about it.

So when I found the paper \cite{nguyen-diaz-nguyen2020}, I decided to type
that embedding theorem and post it on the AMS website of my book, where I put
corrections and related results that I find interesting.

The paper \cite{nguyen-diaz-nguyen2020} led me to another article
\cite{brue-nguyen2020}, and I decided to type that up as well. After that I
discovered the beautiful review paper of Mironescu \cite{mironescu2018}. And
so here we are several hundred pages later.

\smallskip\noindent\textbf{References: }The rule of thumb here is simple: I
only quoted papers and books that I read. I believe that misquoting a paper is
worse than not quoting it. Hence, if a meaningful and relevant article is not
listed in the references, it is because I either forgot to add it or was not
aware of it. While most authors write books because they are experts in a
particular field, I write them because I want to learn a specific topic. I
claim no expertise in fractional Sobolev spaces.

\smallskip\noindent\textbf{Important: }Throughout the book the expression
\begin{equation}
\mathcal{A}\preceq\mathcal{B}\text{\quad means }\mathcal{A}\leq C\mathcal{B}
\label{1d symbol less than}%
\end{equation}
for some constant $C>0$ \textbf{that depends on the parameters quantified in
the statement of the result (usually }$N$, $p$\textbf{ and }$s$\textbf{), but
not on the functions and on their domain of integration}. Also,
\begin{equation}
\mathcal{A}\approx\mathcal{B}\text{\quad means }\mathcal{A}\preceq
\mathcal{B}\text{ and }\mathcal{B}\preceq\mathcal{A}. \label{1d symbol equal}%
\end{equation}
When it is important to stress the explicit dependence of the constant $C$ on
some parameters, I will add a subscript on $\preceq$, for example,
$\preceq_{p,s}$ will mean $\leq C$, where $C=C(p,s)>0$. Finally, when the
constant $C>0$ depends on the domain of integration $I\subseteq\mathbb{R}$ or
$\Omega\subseteq\mathbb{R}^{N}$, I write $\preceq_{I}$ or $\preceq_{\Omega}$, respectively.

\smallskip\noindent\textbf{Web page for appendices, mistakes, comments, and
exercises: }I wrote some appendices for the book \cite{leoni-appendix-FA},
\cite{leoni-appendix-MI}, \cite{leoni-appendix-sobolev}, but, due to the page
limitation (600 pp!), I am posting them on the AMS website for this book.

Also, typos and errors are inevitable in a book of this length and with an
author who is a bit absent minded. I will be very grateful to those readers
who write to giovanni@andrew.cmu.edu to indicate any errors they find. The AMS
is hosting a webpage for this book at

\begin{center}
http://www.ams.org/bookpages/gsm-229/
\end{center}

\noindent where updates, corrections, and other material may be found.

\smallskip\noindent\textbf{Acknowledgments:} I am profoundly indebted to Ian
Tice, who was my go-to person, whenever I got stuck on a proof. Unfortunately
for him, I got stuck a lot! I would also like to thank Kerrek Stinson for
reading parts of the book and for many valuable corrections. Maria Giovanna
Mora and Massimiliano Morini proofread some of the more delicate proofs.
Thanks! I worked with some undergraduate students on parts of this book:
Andrew Chen, Jiewen Hu, Khunpob Sereesuchart, Braden Yates, and Grant Yu. I
hope they had as much fun as I did! I want to thank Xavier Ros-Oton and the
anonymous referees for their feedback and Filippo Cagnetti and Tomasz Tkocz
for useful conversations.

I am grateful to Ina Mette, AMS publisher, Jennifer Wright Sharp, the
production editor, and all the AMS staff I interacted with for their constant
help and technical support during the preparation of this book.

The National Science Foundation partially supported this research under Grants
No.\ DMS-1714098 and 2108784.

Also, many thanks must go to all the people at Carnegie Mellon University's
interlibrary loan for finding the articles I needed quickly. The picture on
the back cover of the book was taken by Kevin Lorenzi.
\aufm{Giovanni Leoni}

\mainmatter

\appendix

\bibliographystyle{amsalpha}
\bibliography{fractional-ref}
\footnote{This is not the complete bibliography. It only containes the articles and books cited in the preface}

%

\end{document}